\documentclass[11pt]{article}
\usepackage{amssymb}
\usepackage{latexsym}
\usepackage{amsmath}

\newtheorem{thm}{Theorem}[section]
\newtheorem{prop}[thm]{Proposition}
\newtheorem{lemma}[thm]{Lemma}

\newtheorem{imprmk}[thm]{Important remark}
\newtheorem{rmk}[thm]{Remark}

\newcommand{\reals}{\mathbb R}

\newcommand{\call}{{\cal L}}

\newcommand{\cinf}{C^{\infty}}

\def\qed{\rule{2.3mm}{2.3mm}}

\newcommand{\lag}{\langle}
\newcommand{\rg}{\rangle}

\begin{document}

\title{\bf Stability of higher order singular
points of Poisson manifolds and Lie algebroids}

\author{
 {\bf Jean-Paul Dufour} \thanks{email: dufourj@math.univ-montp2.fr}\\
 D\'epartement de Math{\'{e}}matiques \\
 CNRS-UMR 5030, Universit\'e Montpellier 2\\
 \\
 {\bf  A\"{\i}ssa Wade} \thanks{email: wade@math.psu.edu}
 \thanks{Research partially supported by the Shapiro Funds} \\
 Department of Mathematics,
  Penn State University \\
 University Park, PA 16802 \\
}
\date{}
\maketitle

\begin{abstract}
We study the stability of singular points
 for smooth Poisson structures  as well as  general Lie algebroids.
 We give sufficient conditions for stability
  lying on the  first-order approximation  (not necessarily linear)
 of the given Poisson structure or Lie algebroid at a singular point.
The main tools used here are the classical Lichnerowicz-Poisson cohomology
 and the deformation cohomology for Lie algebroids  recently introduced
  by Crainic and Moerdijk. We also provide  several examples
 of stable singular points of order $k \geq 1$ for Poisson
 structures and Lie algebroids.
 Finally, we apply our results to
  pre-symplectic leaves of Dirac manifolds. \\

\end{abstract}


\newpage
\section{Introduction and main results}
 Every  bivector field   on a  smooth finite-dimensional
 manifold $M$ will be considered 
 as a  map   $\Lambda: M \rightarrow \bigwedge^2TM$.
  We endow the set ${\cal V}_2(M)$ of
  all smooth bivector fields  with the 
${\cal C}^s$-topology, i.e. the topology of uniform convergence
on compact sets for ${\cal C}^s$-differentiable
  maps and their derivatives up to order $s$. 
 The integer $s$ will be made more precise later on.

A smooth Poisson structure on  $M$ is a
 bivector field $\pi$ that satisfies the equation $[\pi ,\pi ]=0$,
where $[\ ,\ ]$ is the Schouten-Nijenhuis bracket on multi-vector fields
  (see \cite{K85}).
A   {\sl singular point} of $\pi$ is a point $m \in M$ satisfying
 $\pi(m)=0$. Such a point $m$  is said to be {\sl
stable} if every Poisson structure $\widetilde{\pi}$, which is close enough
 to $\pi$, has a singular point near $m.$

 The following result, due to M. Crainic and R.-L. Fernand\`es,
   was announced at  the 4th international conference on Poisson 
 geometry (held in Luxembourg, June 2004):
\begin{thm}{\rm (\cite{CF04})}
Let $\mathfrak {g}$ be the Lie algebra corresponding to the
  1-jet of $\pi$ at a singular point $m_0$.
 If   the second scalar cohomology group of $\mathfrak{g}$
 vanishes then  $m_0$ is stable.
 Conversely, suppose $\mathfrak {g}$ is a Lie algebra  such
  that, every singular point of a degenerate  Poisson
structure whose  1-jet at that point
  corresponds to $\mathfrak {g}$ is stable,
 then  $H^2(\mathfrak{g},\reals)=\{0\}.$  
\end{thm}

After a quick reading of this result, 
one may think that $H^2(\mathfrak{g},\reals)=\{0\}$
  is a necessary and sufficient condition for the stability of a singular
  point of a Poisson structure.
  However, a {\sl curious
phenomenon} occurs: there are  plenty of stable singular points 
  for Poisson structures having a zero 1-jet, hence
 $H^2(\mathfrak{g},\reals )\neq 0$.
In Section \ref{demonstrations}, we will establish a stability
 criterion for general singularities which
improves the one given by Crainic and Fernand\`es.

\medskip
 
 First, we recall some definitions. 
 Two smooth maps  from
$M$ to a manifold $N$ have the same $k$-jet at a point $m\in M$ if they
coincide at $m$ and the derivatives of their local expressions,
 relative to  any local coordinate system,
 agree at that point up to order $k$ ($k\geq 0$).
 We denote by $J^k_{\Phi}(m)$ the {\sl $k$-jet} at $m$ of  a
 function $\Phi \in \cinf(M,N)$.
 Precisely, $J^k_{\Phi}(m)$ is the equivalence class
 of $\Phi$ under the equivalence relation 
 ``to have the same $k$-jet at $m$''.

\medskip
  A {\sl  singularity of order $k$}
 for  a Poisson structure $\pi$ on $M$ is a point $m \in M$, where the
$(k-1)$-jet of $\pi$ vanishes but not its $k$-jet ($k>0$).
 We say that such a singularity  is $k$-{\sl stable} if, 
 for every neighborhood $\Omega$ of $m$ in $M,$ 
there is a neighborhood $\cal W$ of $\pi$ in ${\cal V}_2(M)$
 such that every Poisson structure in $\cal W$ has a
  singularity of order $k$ in $\Omega .$

 Associated to any singular point $m$
 of order $k$ for a Poisson structure $\pi$, there is
  a $k$-homogeneous Poisson structure on $T_mM$ which is
 determined, up to isomorphism, by
 the $k$th-order terms of the Taylor expansion
 of $\pi$ at  $m$. This
 new Poisson tensor $\pi^{(k)}$ will be called
  the {\sl $k$-homogeneous part}  of $\pi$ at $m$.
 We associate to $\pi^{(k)}$  the
 {\sl homogeneous Lichnerowicz-Poisson cohomology complex} (see \cite{DZ04}):
$$\begin{array}{cccccccc}  {\cal V}_1^{(s-k+1)}(T_mM)&&
\stackrel{\partial_1^{^{s-k+1}}}{\longrightarrow} \ \
  {\cal V}_2^{(s)}(T_mM)&& \stackrel
{\partial_2^{^{s}}}{\longrightarrow} \ \
 {\cal V}_3^{(s+k-1)}(T_mM) \cdots \end{array} $$

  \noindent where ${\cal V}_r^{(s)}(T_mM)$ is the space of $s$-homogeneous $r$-vector fields on $T_mM$ (with ${\cal V}_r^{(s)}(T_mM)=\{0\}$ for
  $s<0$),
   and the operators $\partial_r^{\ell}$ are defined by
   $$\partial_r^{\ell}(A)= [\pi^{(k)}, A],$$ for all
  $\ell$-homogeneous multi-vector field $A.$
  The corresponding second cohomology group is
   $$H_{\rm LP}^{2, s}(\pi^{(k)})=
  {{\rm Ker}\left(\partial_2^{^{s}} \right) \over {\rm Im}
   \left( \partial_1^{^{s-k+1}} \right)} \ .$$
   \noindent Under these notations, we have the following theorem,
  where ${\cal V}_2(M)$ is equipped with the ${\cal C}^{2k}$-topology.

\begin{thm} \label{theoreme 1} Let $m$  be a  singularity of order $k$
 for a given Poisson structure whose
  $k$-homogeneous part at $m$ is $\pi^{(k)}$.
   If $H_{\rm LP}^{2, s}(\pi^{(k)})= \{0\}$, 
 for any $s=0, \dots , k-1$, then $m$ is
   a  $k$-stable  singularity. \end{thm}

Our second main result concerns the stability  problem for general
 Lie algebroids. A {\sl Lie algebroid}
 ([CW99]) is a triple $(F, [ \ , \ ],  \varrho)$
 formed by a vector bundle $F$ over a smooth  manifold $M$
  together with a Lie bracket  on the space  $\Gamma (F)$ of smooth sections
 of $F$, and a bundle map $\varrho:  F \rightarrow TM$ called {\sl
 the anchor map} such that
 $$[X, f Y]= f[X, Y]+\varrho(X) \cdot f Y, \quad \forall X, Y \in
 \Gamma (F), \ f \in \cinf(M).$$
The distribution Im$\varrho$ induces a singular foliation on $M$ which is
  called the {\sl foliation of the Lie
algebroid.} A {\sl singular point for a Lie algebroid} is, by definition,
 a point $m$ of the base manifold $M$ where
  the leaf of the associated foliation  reduces to that point.

 Our goal is to give a criterion for the stability of
  these singularities. But first, we have to put a suitable 
 topology on the set of Lie algebroid structures on a given vector bundle.
  To do so and to present our results, we begin
 by recalling that there is a  one-to-one  correspondence  between
 Lie algebroids  and some special Poisson structures.

A  {\sl fiber-wise linear Poisson structure} on a vector bundle $p:E\rightarrow M$ is a Poisson structure on $E$
which has the following properties: \par

\noindent $\bullet$ The bracket of two fiber-wise linear functions is fiber-wise linear.

\noindent $\bullet$  The bracket of a fiber-wise linear function
 and a basic one is basic (every basic function on $E$
  has the form $f\circ p,$  where $f$ is a function on $M$).

\noindent $\bullet$ The bracket of two basic functions vanishes.

 A fiber-wise linear Poisson structure $\Pi$ has the local form
$$\Pi= \sum_{1 \leq i <j \leq r} f_{i,j}^s(x)y_s  \partial y_i \wedge  \partial y_j+
 \sum_{1 \leq i \leq r}\sum_{1 \leq j \leq d } g_{i,j}(x)
  \partial y_i \wedge \partial x_j$$
  in any system of fibered coordinates $(x,y)$, where
 $x=(x_1,\dots ,x_d)$ is a coordinate system
 on the base manifold $M$ and $y=(y_1,\dots ,y_r)$ consists of
 functions that are linear
  on the fibers ($\partial z$ stands for $\partial \over \partial z$).

It is well-known that there is a one-to-one correspondence between Lie
 algebroid structures on a vector bundle $F$
and fiber-wise linear Poisson structures on the total space of its
 dual $E=F^*$ (see for instance [CW99]). Hereafter,
  we identify these two  concepts:
  the concept of a Lie algebroid will  be sometimes employed instead of
 that of a fiber-wise linear Poisson structure.
  In particular, the topology we put on the set of Lie algebroid structures
on a given vector bundle $F$ is the one which corresponds to the 
${\cal C}^s$-topology for bivector fields on the total space of
 the dual vector bundle $E=F^*$
(the integer $s$ will be specified later on).
A point $m$ of the base manifold $M$ is a singular point of a Lie algebroid if and only if $0_m,$ the origin of the
fiber over $m,$ is a singular point of the associated
 fiber-wise linear Poisson structure $\Pi .$ In the above local
expression of $\pi$, such a singular point corresponds to a point
 $x$ with $g_{i,j} (x)=0.$

It makes sense to speak of {\sl  singularities of order $k$  for a Lie 
 algebroid}. By definition, these are singular
points $m$ such that $0_m$ is a  singularity of order $k$  for $\Pi .$ 
The first non-zero terms in the Taylor
expansion of $\Pi$ at $0_m$ have the local form
$$\Pi^{(k)}= \sum f_{i,j}^{(k-1),\ell}(x)y_{\ell}
  \partial y_i \wedge  \partial y_j+
 \sum g_{i,j}^{(k)}(x)
  \partial y_i \wedge \partial x_j$$
where  $g_{i,j}^{(k)}$ (resp. $f_{i,j}^{(k-1), \ell}$) are $k$-homogeneous
 (resp. $(k-1)$-homogeneous) polynomials.
 When $k=1$ we have
$$\Pi^{(1)}= \sum a_{i,j}^{\ell} y_{\ell}
  \partial y_i \wedge  \partial y_j+
 \sum b_{i,j}^{\ell}x_{\ell}
  \partial y_i \wedge \partial x_j$$
  where $a_{i,j}^{\ell}$ are structure constants for a Lie algebra
 $\mathfrak {g}$ and $b_{i,j}^{\ell}$ are constants.
 It corresponds to the algebroid associated to a linear action of
  the Lie algebra  $\mathfrak {g}$ on the tangent space of the base  
 manifold at a given point (see for example [CW99]).

 As for general Poisson structures,
the first homogeneous terms of the Taylor expansion at $0_m$,
  denoted $\Pi^{(k)}$, determine
 a $k$-homogeneous Lie algebroid. 
 Moreover, we can see that $\Pi^{(k)}$ is a
 fiber-wise linear Poisson structure on $V=T_{0_m}E ,$
 viewed as a trivial vector bundle with base $T_mM$ and
 fiber $T_{0_m}(p^{-1}(m))\simeq p^{-1}(m).$
We say that the algebroid singularity $m$ is  $k$-{\sl stable} if, for every neighborhood $\Omega$ of $m$ in $M,$
there is a neighborhood $\cal W$ of $\Pi$ in ${\cal V}_2(E)$ such that every fiber-wise linear Poisson structure in
$\cal W$ has a singularity of order $k$ in $\Omega .$
Even though Lie algebroids can be viewed as
 specific Poisson structures, the Lichnerowicz-Poisson cohomology
  used in Theorem \ref{theoreme 1} is not suitable
  for the study of the stability problem for Lie algebroids.
  One needs to use the {\sl deformation cohomology}   introduced recently
  by Crainic and Moerdijk (see [CM04]).  This is also called {\sl the
  linear fiber-wise Poisson cohomology}
  for Lie algebroids (see [DZ04]). Hereafter,
  we will use  a local and homogeneous version.
We attach to $\Pi^{(k)}$ the  deformation complex
$$\begin{array}{cccccccc}  {\cal V}_{1,lin}^{(s-k+1)}(T_{0_m}E)
&&\stackrel{\partial_{1}^{^{s-k+1}}}{\longrightarrow} \ \
  {\cal V}_{2,lin}^{(s)}(T_{0_m}E)
&& \stackrel{\partial_2^{^{s}}}{\longrightarrow}
 \ \ {\cal V}_{3,lin}^{(s+k-1)}(T_{0_m}E) \cdots \end{array} $$

  \noindent which is by definition the sub-complex of the 
 Lichnerowicz-Poisson complex. Precisely, we replace
 ${\cal V}_r^{(s)}(T_{0_m}E)$
  by its sub-space ${\cal V}_{r,lin}^{(s)}(T_{0_m}E)$ formed by
$s$-homogeneous $r$-vector fields $K$ which satisfies
 the following properties:

\noindent $\bullet$ $K(df_1\dots ,df_r)$  is fiber-wise linear if
 $f_1,\dots ,f_r$, are fiber-wise linear,

\noindent $\bullet$ $K(df_1\dots ,df_r)$  is basic if $f_1,\dots ,f_{r-1}$, 
 are fiber-wise linear and $f_r$ is basic,

\noindent $\bullet$ $K(df_1\dots ,df_r)=0$ if at least two of the
 $f_i$ are basic.

In any local coordinates $(x,y)$ chosen as above, elements of ${\cal V}_{r,lin}^{(s)}(V)$ have the form
$$ \sum f_{i_1\dots i_r}^{(s-1),u}(x)y_u  \partial y_{i_1} \wedge\dots\wedge  \partial y_{i_r}+
 \sum g_{j_1\dots j_{r-1},t}^{(s)}(x)
 \partial x_t\wedge \partial y_{j_1} \wedge\dots\wedge  \partial y_{j_{r-1}} $$
where  $g_{j_1\dots j_{r-1},t}^{(s)}$ are $s$-homogeneous polynomials and
 $f_{i_1\dots i_r}^{(s-1),u}$ are
$(s-1)$-homogeneous polynomials on  $T_mM.$ Denote by
   $H_{lin}^{2, s}(\Pi^{(k)})$ the second group of cohomology for
 the above cohomology complex.

Under these notations we have the following theorem, where  ${\cal V}_2(E)$ is equipped with the ${\cal C}^{2k}$-topology.

\begin{thm} \label{theoreme 2} Let $m$  be a singularity of 
 order $k$  for the algebroid  $\Pi $ and
$\Pi^{(k)}$ its $k$-homogeneous part at $m.$  If $H_{lin}^{2, s}(\Pi^{(k)})= \{0\}$, for any $s=0, \dots , k-1$,
then $m$ is
   a  $k$-stable  algebroid singularity. \end{thm}

Section \ref{demonstrations} contains the proofs of our 
main results (Theorems  \ref{theoreme 1} and \ref{theoreme 2}).
In Section \ref{exemples}, we will give examples and shows why
  most of the singularities of higher order are stable
  (at least in dimension 3). Section \ref {feuilles} gives
  a method for generalizing the above results to singular leaves
 (not necessarily reduced to a point).

\section{Proofs of the  main theorems}\label{demonstrations}

We will start by establishing Theorem \ref{theoreme 1}.
First, one can notice the following

\begin{rmk}
\label{rmk1}
\noindent {\bf a)} The hypothesis of  Theorem \ref{theoreme 1} 
 says that $\partial_2^{s}$ are one-to-one for $s=0, \dots , k-2$, and
$${\rm Ker}\partial_2^{k-1}={\rm Im}\partial_1^{0}.$$

\medskip
\noindent {\bf b)} In contrast with Theorem 1.1, we do not have
  a ``kind of'' converse to Theorem 1.2 when $k>1$. A counter-example is given
 by the Poisson structure
  $\pi^{(2)}= xy {\partial \over \partial x}
 \wedge {\partial \over \partial y} $ 
 on $\reals^2$.
 Every singular point $m$ of order $2$ 
 for a Poisson structure  $\pi$ whose $2$-jet at that point is isomorphic
 to   $\pi^{(2)}$ is 2-stable. Observe that $\pi$ itself is (locally)
 isomorphic to
 $\pi^{(2)}$. However  the Poisson 2-cocycle 
 ${\partial \over \partial x} \wedge {\partial \over \partial y}$ is 
 not a coboundary. Hence, $H_{LP}^{2,0}(\pi^{(2)}) \ne \{0\}$. 
\end{rmk}

 \noindent  Given a neighborhood $\Omega$ of $m$ in $M$,
 we choose a local coordinate system $(x_1,\dots ,x_n),$ 
 defined on an open neighborhood $U$ of $m$
 that is contained $\Omega$.
 For any bivector field $\Lambda$ on $M$ and
 for  any point $p$ in $U$, we denote by $\Lambda^{(s)}_p$ the $s$-homogeneous
 terms of the Taylor expansion of $\Lambda$ at $p$ in the above 
 coordinates. We can view  $\Lambda^{(s)}_p$
 as an element of ${\cal V}_2^{(s)}(\reals^n),$ then the
 {\it $(k-1)$-jet extension} of $\Lambda$ is given by:

 \begin{eqnarray*}
{\cal J}^{k-1}_{\Lambda }:& U  \longrightarrow &
{\cal V}_2^{(0)}(\reals^n)\times\dots\times{\cal V}_2^{(k-1)}(\reals^n)\\
 &p  \longmapsto & \Big(\Lambda^{(0)}_p, \dots , \Lambda^{(k-1)}_p \Big).
\end{eqnarray*}

\noindent Accordingly $\pi^{(k)}=\pi^{(k)}_m$ may be thought of
 as an element of ${\cal V}_2^{(k)}(\reals^n)$ and the associated linear map
$\partial_1^{0}$ sends ${\cal V}_1^{(0)}(\reals^n)$ into
 ${\cal V}_2^{(k-1)}(\reals^n).$
 From now on,  we will use the notation
 $$A=\{0\} \times \cdots \times \{0 \}
  \times {\rm Im}(\partial_1^{0})$$ which is a subspace of
 ${\cal V}_2^{(0)}(\reals^n)\times\dots\times{\cal V}_2^{(k-1)}(\reals^n).$
Then, one gets
 \begin{lemma} \label{Lemme 1} The image of the differential of
 ${\cal J}^{k-1}_\pi $ at $m$  coincides with $A.$
 \end{lemma}

\noindent The proof of this lemma is straightforward.
 It is left to the reader.

\medskip

\noindent
 Given a point $p$ in $U$ and bivector field $\Lambda$ on $M,$ we define
 $${\cal F}_{p,\Lambda}:{\cal V}_2^{(0)}(\reals^n)\times\dots\times{\cal V}_2^{(k-1)}(\reals^n)\longrightarrow
 {\cal V}_3^{(k-1)}(\reals^n)\times\dots\times{\cal V}_3^{(2k-2)}(\reals^n)$$
 such that the components of ${\cal F}_{p, \Lambda}( v_0 , \cdots , v_{k-1})$ are given by
   $${\cal F}_{p, \Lambda}^{\ell}( v_0 , \cdots , v_{k-1}) =
   \sum_{i\leq \ell-k} [v_ {i}, \Lambda^{(\ell-i)}_p]+
   \ { 1 \over 2} \sum_{\ell-k < i,j \leq \ell, \ i+j=\ell} [ v_ i , v_j],$$
   \noindent for all $k \leq \ell \leq 2k-1,$ and where 
$[\ ,\ ]$ is the Schouten-Nijenhuis bracket. One has:
   \begin{lemma} \label{Lemme 2} 
 Under the hypothesis of Theorem \ref{theoreme 1},
 $A$ is the kernel of the differential
 of ${\cal F}_{m, \pi}$ at the origin $\mathbf{O}$.
 \end{lemma}
\noindent {\it Proof:}   When $k=1$, one gets 
${\cal F}_{m, \pi}( v_0)= [v_0, \pi^{(1)}].$
Differentiating, one can easily check that
${\rm Ker}\big( ({\cal F}_{m, \pi})_* (\mathbf{O})\big)
={\rm Ker}(\partial_2^{^0})
={\rm Im}(\partial_1^{^0}).$
 Now, assume that $k>1$, then the equation
$$\big(({\cal F}_{m, \pi})_*(\mathbf{O}\big)\big)(w_0, \dots , w_{k-1})=0$$ 
  gives the system of equations
$$
[w_{0}, \pi^{(k)}]=0, \quad
 \quad [w_1 , \pi^{(k)}]+ [w_0, \pi^{(k+1)}]=0, \ \dots,$$
$$[w_{k-1}, \pi^{(k)}]+ \cdots + [w_{0}, \pi^{(2k-1)}]=0
$$
\noindent The first equation
 of this system  says that $w_0 \in {\rm ker}(\partial_2^{0})$.
By Remark \ref{rmk1},  one gets $w_0 \in {\rm Im}(\partial_1^{1-k}).$
 Using the fact that  ${\rm Im}(\partial_1^{^{1-k}})=\{0\}$ for $k>1$,
 one obtains $w_0=0$. By an iteration procedure,
 one gets
$$w_0= \cdots = w_{k-2}=0 \quad {\rm and} \quad \partial_2^{^{k-1}}(w_{k-1})=0,
 \quad  \forall \ k>1.$$ 
\noindent This shows that
 the kernel of the
 differential of ${\cal F}_{m, \pi}$ at $\mathbf{O}$ is exactly $A$. 

\hfill \qed

   \noindent {\bf Proof of Theorem \ref{theoreme 1}:}\par

\noindent Once for all, we fix a complement subspace $B$ of
   $A,$ i.e.
  $${\cal V}_2^{(0)}(\reals^n)\times\dots\times{\cal V}_2^{(k-1)}(\reals^n) =A\oplus B.$$

  \noindent Lemma \ref{Lemme 2} implies that the differential
   of ${\cal F}_{m, \pi}$ at the origin is one-to-one, when restricted to $B$.
 Applying the implicit function theorem,   one gets
   that ${\cal F}_{m, \pi}$  is an injective immersion when it is restricted
  to a neighborhood ${\cal B}_0$ of the origin
 in $B$. Note that ${\cal F}_{p, \Lambda}$ depends continuously on $\Lambda$
 in the ${\cal C}^{2k}$
 topology. A classical singularity technique (see Lemma A, p. 61, 
 \cite{GG73}) implies that there is a
   neighborhood ${\cal B}_1\subset{\cal B}_0$ of the origin in $B$,
   a neighborhood $\Omega_1 \subset U$ of $m$,  and a neighborhood
   ${\cal W}_1$ of $\pi$ in ${\cal V}_2(M)$ such that, for every $p$
 in $\Omega_1$ and for every $\Lambda$ in ${\cal W}_1,$
${\cal F}_{p,\Lambda}$ is one-to-one when restricted to ${\cal B}_1.$

It follows from Lemma \ref{Lemme 1}  that the map ${\cal J}^{k-1}_\pi $ is
 transversal to $B$ at the origin. 
By transversality, there is a neighborhood
 ${\cal W}_2$ of $\pi$ in ${\cal V}_2(M)$
such that, for every $\Lambda$ in ${\cal W}_2,$ ${\cal J}^{k-1}_\Lambda $
 intersects (transversally) ${\cal B}_1$
at a point $p\in \Omega_1 .$ 
 Set ${\cal W}={\cal W}_1\cap {\cal W}_2.$
Now pick a {\em Poisson structure} $\Lambda$ in ${\cal W}.$ We know that
 there is a point $p$ in $\Omega_1$ where ${\cal J}^{k-1}_\Lambda (p)$
  intersects ${\cal B}_1.$ Considering the
equation $[\Lambda ,\Lambda ]=0$, we can notice  that the terms of order
 $k-1,\dots ,2k-2$ in the Taylor expansion
of  $[\Lambda ,\Lambda ]$ at that point $p$ give
$${\cal F}_{p, \Lambda}({\cal J}^{k-1}_\Lambda (p))=0.$$
The injectivity of $({\cal F}_{p, \Lambda})|_{{\cal B}_1}$
 implies $${\cal J}^{k-1}_{\Lambda }(p)=0.$$
If necessary, we may replace $\cal W$ by  a smaller neighborhood of $\pi$
  to ensure that the $k$-jet of $\Lambda$ at $p$ is non-trivial.
 So, $p$ is the singularity point for $\Lambda$ we sought. 
 This completes the proof of  Theorem \ref{theoreme 1}.

\hfill \qed

\begin {imprmk} 
 In our study there is no need for the singularity to be isolated. In fact,
 from the proof of Theorem \ref{theoreme 1}, one sees that
   the set of singularities  of order $k$  for
  Poisson structures $\Lambda$
  sufficiently close to $\pi$ is (locally) formed  by points $p$ such that 
${\cal J}^{k-1}_{\Lambda }(p) \in {\cal B}_1 $, for each $k$  fixed.
  So, this set is, by transversality,
 a sub-manifold of dimension $r=dim(A)$. One
can also notice that this manifold depends continuously on $\Lambda $,
 in a natural sense.
\end{imprmk}

   The proof of Theorem \ref{theoreme 2} is a direct adaptation
 of that of Theorem \ref{theoreme 1}.
Essentially, one replaces the symbol ${\cal V}_u^{(s)}$ by
  ${\cal V}_{u,lin}^{(s)}$ and
$${\cal J}^{k-1}_{\Lambda } : U\longrightarrow
 {\cal V}_2^{(0)}(\reals^n)\times\dots\times{\cal V}_2^{(k-1)}(\reals^n)$$
 $$p\longmapsto (\Lambda^{(0)}_p, \dots , \Lambda^{(k-1)}_p)$$  by
$${\cal J}^{k-1}_{\Lambda}: U \longrightarrow
 {\cal V}_{2,lin}^{(0)}(\reals^d\times\reals^r)\times\dots\times{\cal V}_{2,lin}^{(k-1)}(\reals^d\times\reals^r).$$
$$ p\longmapsto (\Lambda^{(0)}_{0_p}, \dots , \Lambda^{(k-1)}_{0_p}),  \ \ \ \ 
\ \ \ \ \ \ \ \ \ \ \ \ \ $$
\noindent where $d$ and $r$ are the dimensions
 of the base manifold and the fiber, respectively.

\section{Examples}\label{exemples}

\subsection{Singularities of Poisson structures in dimension 2}

In dimension 2, singularities of Poisson structures are singularities
 of functions for which  there are specific classical techniques
 (see \cite{M01}). For instance,
  singularities of topologically stable Poisson structures in  Radko's
sense (see \cite{R02}) are 1-stable.

\subsection{Singularities of order 1 for Poisson structures}

Let $m$ be a singularity of  order 1 for
  a Poisson structure $\pi .$ Denote by
 ${\mathfrak g}$ the Lie algebra corresponding to the
linear part of $\pi$ at $m.$ In this case,
  the hypothesis of Theorem
 \ref{theoreme 1} reduces to $H_{\rm LP}^{2, 0}(\pi^{(1)})= \{0\}.$ 
 A straightforward calculation shows that this  condition is equivalent
 to $H^2({\mathfrak g},\reals )=0.$  This is exactly the condition
 given by Crainic and Fernand\`es (see Section 1).

On can notice that  if the Lie algebra associated to the
linear part of $\pi$ at $m$ is strongly rigid  in the sense of
 Bordemann, Makhlouf, and  Petit (see \cite{BMP02})
 then $m$  is 1-stable.

\subsection{Singularities of order 2 for Poisson structures}

It is shown in [DH91] that  Poisson structures on a vector space 
$V$  having the  form
$$\pi^{(2)}=\sum a_{i,j}x_ix_j\partial x_i\wedge\partial x_j$$
\noindent are generic among quadratic Poisson structures on $V$.
We consider such a Poisson structure $\pi^{(2)}$. It is known that the
 $a_{i,j}$ are invariants for this Poisson structure (see \cite{DW98}).
 Set $\lambda_i=\sum_{j}a_{i,j}.$

\begin{lemma}  If $\lambda_i\ne 0$ for every $i$,
  and $\lambda_i+\lambda_j\ne 0$ for  $i<j$.
   then $H_{\rm LP}^{2, s}(\pi^{(2)})= \{0\}$, for $s=0$ and $s=1.$
\end{lemma}
Note that the hypothesis in the above lemma
  may occurs only  when the dimension is at least 3.
The proof of this lemma can be found in
 \cite{M01}. This lemma and Theorem \ref{theoreme 1}
 imply that that singularities  of order 2  of
  quadratic Poisson structures  are,
in general, stable (in dimension $\geq 3$).

\subsection{Singularities of Poisson structures in dimension 3}

Let $\omega$ be a volume form on a 3-dimensional manifold $M$.
 The map $\pi\mapsto i_\pi\omega$
establishes a one-to-one correspondence between Poisson structures
 and integrable 1-forms (i.e. 1-forms $\alpha$ such
that $\alpha\wedge d\alpha =0$). Moreover, $\pi$ and $\alpha =i_\pi\omega$
  share the same (singular) foliation.
 In \cite{CN82},  Camacho and  Lins-Neto 
  proved that singularities of order $k$  for these integrable forms
 are, in general, stable. Hence, the same conclusion
 holds for Poisson structures in dimension 3.

  In fact, the proof of Theorem \ref{theoreme 1}  
 is inspired by  methods used for integrable 1-forms. 
 We  will interpret, in terms of Poisson structures,
 two results for integrable 1-forms  obtained by Camacho and  Lins-Neto.   
First, recall that a   vector field 
 $X=X_1\partial x+X_2\partial y+X_3\partial z$  on $\reals^3$, which vanishes
at the origin, has an {\sl algebraically isolated zero} 
 at the origin if the ideal generated by the
germs of its components $X_1,$ $X_2$ and $X_3$
 has a finite codimension in the space of germs of functions of $x,$
$y$ and $z$ at the origin.

\begin{lemma} \label{lemmeLN} Let $\pi^{(k)}$ a $k$-homogeneous Poisson structure on $\reals^3.$
 Suppose its modular
vector field $X^{(k-1)}$ with respect to $\omega=dx\wedge  dy\wedge dz$
 has an algebraically isolated zero at the
origin. Then the hypothesis of Theorem \ref{theoreme 1} 
 (i.e $H_{\rm LP}^{2, s}(\pi^{(k)})= \{0\}$, for any $s=0,
\dots , k-1$) holds true.\end{lemma}

Hence any order $k$ singularity, with $\pi^{(k)}$ as in this lemma,
 is $k$-stable.

\begin{lemma} Under the hypothesis  of Lemma \ref{lemmeLN} and
 for all $k>2$, we have
$$\pi^{(k)}= \frac{1}{k+1}I\wedge X^{(k-1)},$$
where $I=x\partial x+y\partial y+z\partial z$ is the Liouville vector field.
 
\end{lemma}

This last lemma allows to construct many examples of
 $k$-stable singularities for  $k>2$. For instance,
the origin is a $k$-stable singularity for
$$\pi =(xy^{k-1}-yx^{k-1})\partial x\wedge \partial y+
(yz^{k-1}-zy^{k-1})\partial y\wedge \partial z+
(zx^{k-1}-xz^{k-1})\partial z\wedge \partial x.$$

\subsection{Singularities of order 1 for Lie algebroids}

We recall that if $m$ is a singularity order 1
 for a Lie algebroid then the linear part of this algebroid at that
singularity is a Lie  algebroid corresponding to a linear action of a
  Lie algebra $\mathfrak{g}$ on a vector space
$V.$ We will see that, in this case,  the condition of stability of
 Theorem \ref{theoreme 2} is equivalent to $H^1(\mathfrak{g},V^*)=0.$

\medskip
 Consider a fiber-wise linear Poisson structure that
 corresponds to an action Lie  algebroid
 $g \ltimes V$. Suppose that it is  given
 by the following expression in a fibered coordinate system 
$(x_1, \dots , x_d,y_1, \dots , y_r)$:
 
 $$\Pi^{(1)}= \sum_k a_{i,j}^ky_k  \partial y_i \wedge  \partial y_j+
 \sum_{\ell} b_{p ,q}^{\ell}x_{\ell}  \partial y_p \wedge \partial x_q,$$
  where $a_{i,j}^k$ are structure constants of the Lie algebra
 $\mathfrak {g}$ and $b_{p,q}^{\ell}$ are constants which determine the action
  $\mathfrak {g}$ on $V$.
 Set $$\mu= \sum \mu_{u,v} \partial y_u \wedge \partial x_v,$$
 \noindent where the $\mu_{u,v}$ are constants. Consider the equation
 $[\Pi^{(1)}, \mu]=0$. In this equation, the coefficient
 of $ \partial x_q \wedge \partial y_i \wedge \partial y_j$ gives
$$ \sum_v \Big(\mu_{i,v}b_{j,q}^v -  \mu_{j,v}b_{i,q}^v \Big)-
 \sum_u a_{i,j}^u \mu_{u,q}=0.
 \leqno{(1)}$$
\noindent But, if $V^*$
 is identified with the vector space of constant
 vector fields generated by the vector fields $\partial x_i$, then
 the dual action 
 $\varrho : \mathfrak{g} \rightarrow {\rm End}(V^*)$ is
 given by  $$\varrho(y_i)(\partial x_j)= -
 \sum_v b_{i, v}^j \partial x_v.$$
\noindent Furthermore, if we set
 $$\mu_i= \sum_u \mu_{i, \ell} \partial x_{\ell},$$
\noindent  we obtain that Equation (1) can be expressed as
 $$ \langle \varrho(y_i) (\mu_j)- \varrho(y_j) (\mu_i),  \ dx_q \rangle
 = \sum_u a_{i,j}^u \mu_{u,q}.$$
\noindent More precisely, Equation (1) says that
  $\mu$ corresponds to a 1-cocycle for the representation $\varrho$.
 This shows that constant 2-cocycles for the deformation cohomology
 complex can be identified with 1-cocycles relative
 to the representation $\varrho$.
 Similarly, 2-coboundaries for the deformation cohomology complex
 can be identified with 1-coboundaries 
  relative to the representation $\varrho$.
 There follows that $$H_{lin}^{2, 0}(\Pi^{(1)}) \cong H^1(\mathfrak{g},V^*).$$
In particular, when $\Pi^{(1)}$ is the  fiber-wise Poisson
structure induced by the natural action of $\mathfrak{g}$ on itself, then
 $H_{lin}^{2, 0}(\Pi^{(1)}) \cong H^1(\mathfrak{g},\mathfrak{g}^*)$.
 As R.-L. Fernand\`es points out (private communication),
    one  can notice that, in this case,
  the hypothesis of Theorem \ref{theoreme 2}  is stronger than  
 that of Theorem \ref{theoreme 1}. More precisely,  the  natural inclusion 
 $H^2(\mathfrak{g}, \reals)  \hookrightarrow H^1(\mathfrak{g},\mathfrak{g}^*)$
 is not onto, in general 
 (since $H^2(\mathfrak{g}, \reals)$  can be viewed as the skew-symmetric
 part of  $H^1(\mathfrak{g}, \mathfrak{g}^*$)).
 Consequently, the fact that $m$ is a 1-stable singular point for a Poisson 
 structure whose 1-jet at $m$ corresponds to $\mathfrak{g}$ does not
  guaranty that it is 1-stable for the associated Lie algebroid.
For instance,   if $ \mathfrak{g}= \mathfrak{aff}(1)$  (i.e. the 
 Lie algebra of affine transformations on $\reals^2$) then the origin
 $O$ is 1-stable for the associated Lie-Poisson structure on 
 $\mathfrak{g}^*$ since $H^2(\mathfrak{aff}(1), \reals)= \{0\}$.
 But, the origin is not 1-stable  for the corresponding Lie algebroid.

\bigskip

Furthermore, one can notice that, in the particular case where 
$\mathfrak{g}$ is reductive and  the $X_i=  \sum b_{ik}^{\ell}
 x_{\ell} \partial x_k$ are diagonal, the 
Hochschild-Serre factorization theorem (see \cite{HS53}) gives
$H_{lin}^{2,0}(\Pi^{(1)})=\{0\} \iff
 H^1(\mathfrak g, \mathbb{R})=\{0\}.$

\subsection{Singularity of order 2 for Lie algebroids}
Consider $\reals^{n+1}$ ($n>1$)  with the standard coordinates
 $(x_0, x_1, \dots , x_n)$ and the function
 $\rho= x_0^2+ \dots + x_n^2$. Let $\Pi$ be the fiber-wise linear
  Poisson structure of on $\reals^{2n+2}$ given by
 $$\Pi^{(2)}= (\rho-1) \sum_{i=0}^n \partial x_i   \wedge \partial y_i
 +  \sum_{i<j} (x_i y_j - x_j y_i) \partial y_i   \wedge \partial y_j.$$
\noindent  This is the Lie algebroid
 attached to $(\rho-1)$, in the sense of Monnier (see [M02]).
 Every point lying on  the  unit sphere  is singular. 
 It follows  from results proven in \cite{M01}
  that $H_{\rm lin}^{2,s}( \Pi^{(2)})$ vanishes for all
 $s$ (see Lemma 5.4.3 and Proposition 6.1.2 in  \cite{M01}).
  Therefore, the unit sphere is 2-stable.

\medskip

\section{Application: stability of non-trivial leaves}\label{feuilles}

In this section, we  will apply our previous results to Dirac manifolds
 and non-trivial leaves of Lie algebroids.
Let $M$ be a smooth $n$-dimensional manifold. We consider the vector bundle
 $TM \oplus T^*M$ together with the Courant bracket defined by 
$$[(X_1, \alpha_1), (X_2, \alpha_2)]_{_C}=([X_1,X_2],
  \call_{X_1}\alpha_2 - i_{X_2}d\alpha_1),$$
for any  $(X_1, \alpha_1), \ (X_2, \alpha_2) \in \Gamma(TM \oplus T^*M)$.
 Recall that a {\sl Dirac structure} on $M$ 
 (see \cite{C90}) is a sub-bundle $L \subset TM \oplus T^*M$  of rank $n$
 which is closed under the Courant bracket and which is
isotropic with respect to the bilinear symmetric operation given by
$$ \lag (X_1, \alpha_1), (X_2, \alpha_2) \rg ={1 \over 2}(i_{X_2} \alpha_1 + i_{X_1} \alpha_2),$$
 for any $(X_1, \alpha_1),  \ (X_2, \alpha_2) \in
 \Gamma(TM \oplus T^*M).$ In this case, $(M, L)$ is called
 a {\sl Dirac manifold}. It is known that every Dirac manifold admits
 a  foliation  by pre-symplectic leaves (see \cite{C90}).
 We also know that the dimensions of the leaves have the same parity 
 (see \cite{DW04}).
 
\medskip
 A  pre-symplectic leaf $S$ of $(M,L)$ is {\sl singular}
  if every  neighborhood of $S$ intersects leaves of higher dimension.
 We say that  the singular leaf 
$S$ is {\sl stable} if, for any  neighborhood $U$ of 
 $S$,  there is a neighborhood  ${\cal W}$ of $L$ such that  all 
 Dirac structures in  ${\cal W}$  admit a 
 pre-symplectic leaf near $S$ having the same dimension as $S$.
 Here, we endow the set of Dirac structures with the $C^s$-topology by means of the $C^s$-topology on local smooth sections
 of $TM \oplus T^*M$ ($s$  will be specified below).
 
 Let $S$  be a pre-symplectic leaf of $(M,L)$, it is proven in 
 \cite{DW04} that if $N\subset M$  is an embedded submanifold 
 that intersects transversally 
 $S$ at $m_0$ (in the sense that $T_{m_0}M=T_{m_0}S \oplus T_{m_0}N$)
 then $L$ induces
  a Poisson structure on $N$ which vanishes at $m_0$. Moreover,
 this transverse Poisson structure does not depend on $m_0$:
 it is unique up to Poisson isomorphisms.
 This is called {\sl the transverse Poisson structure along $S$}.
 In the following theorem, we assume that $m_0$ is 
 $k$-stable and we endow the set of Dirac structures with the
 $C^{2k}$-topology.

\begin{prop} Assume that $S$ is an embedded singular  pre-symplectic leaf
 of a Dirac manifold $(M,L)$ and $m_0 \in S$ is a $k$-stable singular point
   for the transverse  Poisson structure. Then $S$ is stable.
\end{prop}

\noindent{\it Proof:} Consider a  tubular neighborhood
 of $S$ that corresponds to the vector bundle $p : E \rightarrow S$.
 Denote by  $\pi_{_V}$ the transverse Poisson structure defined on
 an open neighborhood $U \subset E_{m_0}$ of $m_0$
  with $\pi_{_V}(m_0)=0$. We set $N=E_{m_0}$. Any Dirac structure
 $L'$   sufficiently close to $L_{|E}$ has  pre-symplectic
 leaves $S'$  near $S$ that intersect  $N$ with the transversality relation:
  $$T_{m}E= T_{m}S' + T_{m} N $$ at  points $m \in U \subset N$.
  Using the transversality relation and methods developed in \cite{DW04},
 one can show that $L'$ induces a  Dirac structure
 $L'_{_V} $ in a neighborhood of $m_0$ in $N$.
  Moreover, the fact that $L'$ is chosen 
 to be very close to $L_{|E}$
 implies that  $L'_{_V} $ is  the graph of a Poisson structure
$\pi'_{_V}$ which is close to  $\pi_{_V}$ on $U$
  (if necessary, we replace $U$ by a smaller open set).
 Since $m_0$ is $k$-stable for  $\pi_{_V}$,
  the Poisson structure $\pi'_{_V}$ has a singular point 
  of order $k$ in $U$. Obviously, singular points
  for  $\pi'_{_V}$ are located on leaves of the Dirac manifold
 $(E, L')$ that have the same dimension as  $S$. In other words,
 if  $m_1 \in S' \cap U$ is   a singular point for  $\pi'_{_V}$ then
 we have a direct sum
 $$T_{m_1}E= T_{m_1}S' \oplus T_{m_1} N,$$ moreover  $S'$  projects
 diffeomorphically onto $S$ when restricted in a small neighborhood of $m_0$
  in $M$. 
 
\hfill \qed

\medskip

\begin{rmk}  In the above proof,  
 $S'$ does not globally project onto $S$ via the tubular neighborhood
 projection, in general (specially, when $S$ is non-compact). 
  For instance, let
  $M$ be the product of  $\reals^2$
 by the torus $\mathbb T^2$ and suppose that
  $L_{\pi}$ is the graph of the Poisson structure ${\pi}$ given by 
 $$\pi= y {\partial \over \partial x} \wedge 
{\partial \over \partial y} + {\partial \over \partial \theta_1} \wedge 
{\partial \over \partial \theta_2}.$$ The submanifolds given by equations
   $x=$cte, $y=0$ are singular leaves of $(M, L_{\pi})$.
 By Proposition 4.1, each of these leaves is stable.
 Note that they are diffeomorphic to  $\mathbb T^2$.
Now, we consider 
$$\pi_{\varepsilon}= y {\partial \over \partial x} \wedge 
{\partial \over \partial y} + 
\Big( {\partial \over \partial \theta_1}  + \varepsilon 
 {\partial \over \partial x} \Big) \wedge 
{\partial \over \partial \theta_2},$$
where $\varepsilon \in \reals $ is a deformation parameter.
 For $\varepsilon \ne 0$, every
 2-dimensional  symplectic leaf of $\pi_{\varepsilon}$
 is (globally) diffeomorphic to $S^1 \times \reals$.
  Hence, it is not (globally) diffeomorphic to the torus. 

\medskip
We should mention that
  Crainic and Fernand\`es have obtained (private communication)
 results on 1-stability of compact leaves of Poisson structures 
 in a stronger sense, i.e.  $S$  diffeomorphic to $S'$.

\hfill \qed
\end{rmk}

Using the notion of transverse Lie algebroid structure along a leaf 
 (see, for instance, \cite{F02}) and the above techniques, one can similarly
 show the following proposition:

\begin{prop} Let $S$ an embedded singular leaf of an algebroid structure.
 If $m$ is a k-stable singularity of
the transversal algebroid structure along  $S$,
 then $S$ is stable.\end{prop}

\noindent{\bf Aknowledgement.} We would like to thank
 Rui  Loja Fernand\`es for  fruitful
 discussions and  comments on an earlier
 version of this paper.

\end{document}